\documentclass[11pt]{article}
\usepackage[T1]{fontenc}
\usepackage{lmodern,microtype}
\usepackage{amsmath,amssymb,amsthm,mathtools,mathrsfs}
\usepackage[margin=1in]{geometry}
\usepackage{enumitem}
\usepackage[colorlinks=true,citecolor=blue,linkcolor=blue,urlcolor=blue]{hyperref}
\setlist[enumerate]{label=\textup{(\roman*)},leftmargin=*}

\newtheorem{maintheorem}{Theorem}

\newtheorem{maincorollary}[maintheorem]{Corollary}
\newtheorem{theorem}{Theorem}[section]
\newtheorem{proposition}[theorem]{Proposition}
\newtheorem{lemma}[theorem]{Lemma}
\newtheorem{corollary}[theorem]{Corollary}
\newtheorem*{introconjecture}{Conjecture}
\theoremstyle{remark}
\newtheorem{remark}[theorem]{Remark}
\newtheorem*{introremark}{Remark}
\newcommand{\papersection}[1]{\section{#1}}
\usepackage{titlesec}

\titleformat{\section}[block]
  {\normalfont\Large\bfseries\centering}
  {\thesection.}
  {0.5em}
  {}
\titleformat{\subsection}
  {\normalfont\large\bfseries\raggedright}
  {\thesubsection.}
  {0.5em}
  {}

\newcommand{\Cbb}{\mathbb C}
\newcommand{\Qbb}{\mathbb Q}
\newcommand{\Ball}{\mathbb B}
\newcommand{\Xbar}{\overline{X}}
\newcommand{\Ybar}{\overline{Y}}
\newcommand{\Pcal}{\mathcal P}
\newcommand{\Lcal}{\mathcal L}
\newcommand{\Ocal}{\mathcal O}
\newcommand{\CC}{\operatorname{CC}}
\newcommand{\Perv}{\operatorname{Perv}}
\newcommand{\Supp}{\operatorname{Supp}}
\newcommand{\Ilog}{I^{\mathrm{log}}}
\newcommand{\sys}{\operatorname{sys}}

\newcommand{\PP}{\mathbb P}

\title{Euler Characteristics of Perverse Sheaves\\on Complex Ball Quotients}
\author{Mert Akdenizli and Soheil Memariansorkhabi}
\date{}
\begin{document}
\maketitle
\begin{abstract}
Let $X=\Gamma\backslash\Ball^n$ be a finite-volume complex ball quotient
with neat lattice. We prove that every nonzero algebraically constructible
perverse sheaf on $X$ has positive Euler characteristic. More generally,
we establish a criterion for strict positivity on smooth quasi-projective
varieties whose logarithmic cotangent bundle is nef and admits an ample
negative rational twist of the boundary divisor.

For perverse sheaves with positive-dimensional support, we give a lower
bound in terms of the intrinsic geometry of $X$ for both proper and
nonproper quotients. The bound depends on the systole of $X$ and tends
to infinity with it. In particular, the minimum Euler characteristic
among perverse sheaves with positive-dimensional support tends to infinity
as one passes to higher levels of any cofinal normal tower of finite covers.
We also deduce lower bounds for the middle intersection Betti numbers of
projective subvarieties and boundedness of supports with bounded middle
intersection cohomology.
\end{abstract}

\pdfbookmark[1]{Introduction}{introduction}
\section*{Introduction}

Hopf's work in 1931 and Singer's in the 1970s led to the following
conjecture on the Euler characteristic of closed aspherical manifolds: 

\begin{introconjecture}[Singer--Hopf]
Let $M$ be a closed aspherical manifold of real dimension $2n$.
Then
\[
  (-1)^n\chi_{\mathrm{top}}(M)\geq0.
\]
\end{introconjecture}

More generally, Liu--Maxim--Wang
\cite[Conjecture~6.2]{LiuMaximWang} proposed a
perverse-sheaf version of this conjecture for aspherical
smooth projective varieties, subsequently it is extended to
compact aspherical K\"ahler manifolds by Arapura--Wang
\cite[Conjecture~1.2]{ArapuraWang}.
The conjecture asserts that every
$\Pcal\in\Perv(X,\Cbb)$ on such a manifold satisfies
$\chi(X,\Pcal)\geq0$, where
\[
  \chi(X,\Pcal)
  :=\sum_k(-1)^k\dim_{\Cbb}\mathbb H^k(X,\Pcal).
\]
The Euler characteristic of a perverse sheaf on a smooth
projective variety is nonnegative when the cotangent bundle
is nef, and strictly positive when the cotangent bundle is
ample and the sheaf is nonzero \cite{LiuMaximWang,ArapuraWang}.

Let $X=\Gamma\backslash\Ball^n$ be a finite-volume
complex ball quotient with neat lattice $\Gamma$.
In the compact case, the cotangent bundle of $X$ is ample,
so the preceding result gives
$\chi(X,\Pcal)>0$ for every nonzero
$\Pcal\in\Perv(X,\Cbb)$.
We establish the same strict positivity statement
in the noncompact case:

\begin{maintheorem}\label{thm:positive}
Let $X=\Gamma\backslash\Ball^n$ be a noncompact
finite-volume complex ball quotient by a neat lattice.
Then
\[
  \chi(X,\Pcal)>0
\]
for every nonzero $\Pcal\in\Perv(X,\Cbb)$.
\end{maintheorem}

Theorem~\ref{thm:positive} is a consequence of the following general
criterion. It applies to any smooth quasi-projective variety with a
smooth projective SNC compactification.

\begin{maintheorem}
\label{thm:criterion}
Let $\Ybar$ be a smooth complex projective variety and let
$D=\sum_{i=1}^{r}D_i$ be a reduced simple normal crossings divisor.
Set $Y=\Ybar\setminus D$ and $E=\Omega^1_{\Ybar}(\log D)$.
Suppose $E$ is nef and there exist positive rational numbers
$a_1,\ldots,a_r$ such that the $\Qbb$-twisted vector bundle
\(
 E\left\langle-\sum_{i=1}^{r}a_iD_i\right\rangle
\)
is ample. Then every nonzero $\Pcal\in\Perv(Y,\Cbb)$ satisfies
\[
 \chi(Y,\Pcal)>0.
\]
\end{maintheorem}

Let $\Xbar$ be the smooth projective toroidal compactification of $X=\Gamma\backslash\Ball^n$ and
$D=\Xbar\setminus X$ its boundary divisor. For $n\geq2$, neither
the cotangent bundle $\Omega^1_{\Xbar}$ nor the logarithmic
cotangent bundle $\Omega^1_{\Xbar}(\log D)$ is ample. Indeed,
each boundary component is an abelian variety, and the restrictions
of both bundles to it have nonzero trivial quotients.
Nevertheless, $E=\Omega^1_{\Xbar}(\log D)$ is nef by
\cite[Corollary~G]{Memarian}, and $E\langle-aD\rangle$ is ample
for suitable positive rational $a$ by \cite[Theorem~A]{Memarian}.
Theorem~\ref{thm:criterion} therefore proves
Theorem~\ref{thm:positive}. This gives strict positivity of the perverse Euler
characteristics on every finite-volume complex ball quotient with
neat lattice.

Moreover, we obtain an explicit lower bound for the Euler
characteristic of a perverse sheaf in terms of the intrinsic
geometry of $X$. Let $\sys(X)$ denote the systole, the length of the shortest
closed geodesic with respect to the Bergman metric normalized
to have holomorphic sectional curvature $-1$; see Section~\ref{sec:complete}. Our main 
result is the following:
\begin{maintheorem}\label{thm:systolic}
Let $X=\Gamma\backslash\Ball^n$ be a finite-volume
complex ball quotient with neat lattice, and let
$\Pcal\in\Perv(X,\Cbb)$ have an irreducible component
of support of dimension $m>0$.
\begin{enumerate}
\item If $X$ is proper, then
\[
  \chi(X,\Pcal)\geq
  \sinh^{2m}\!\left(\tfrac14\sys(X)\right).
\]
\item If $X$ is nonproper and
\(
  \sys(X)\geq
  4\log\bigl(5n+(4\pi)^4\bigr),
\)
then
\[
  \chi(X,\Pcal)>
  \left(\frac{e^{\sys(X)/16}}{4\pi}\right)^m.
\]
\end{enumerate}
\end{maintheorem}

Let $Z\subset X$ be an irreducible subvariety of dimension $m>0$,
possibly singular or noncompact. We write $IH^q(Z,\Qbb)$ for its
intersection cohomology group in degree $q$ with rational
coefficients. These groups agree with ordinary cohomology when
$Z$ is smooth \cite[\S6.7]{MaximIntersectionHomology}. Let $i:Z\hookrightarrow X$ be the inclusion and let $IC_Z$ be
the intersection complex of $Z$, normalized to be perverse.
Then
\[
  \chi(X,i_*IC_Z)=(-1)^m\chi_{\mathrm{IH}}(Z),
  \qquad
  \chi_{\mathrm{IH}}(Z)
  =\sum_q(-1)^q\dim_{\Qbb}IH^q(Z,\Qbb).
\]
Theorems~\ref{thm:positive} and~\ref{thm:systolic} therefore give
estimates for the signed intersection cohomology Euler
characteristic, including when $Z$ is noncompact. If $Z$ is
smooth, then $IC_Z=\Cbb_Z[m]$ and
$\chi_{\mathrm{IH}}(Z)=\chi_{\mathrm{top}}(Z)$.
The dimensions
\[
  b_q^{\mathrm{IH}}(Z)=\dim_{\Qbb}IH^q(Z,\Qbb)
\]
are called the intersection Betti numbers. They agree with the
ordinary Betti numbers when $Z$ is smooth. For projective
subvarieties, hard Lefschetz gives the following bounds on
middle cohomology:

\begin{samepage}
\begin{maincorollary}[Bounds for middle Betti numbers]
\label{cor:betti}
Let $X=\Gamma\backslash\Ball^n$ be a finite-volume complex ball
quotient with neat lattice, and set $s=\sys(X)$.
If $X$ is noncompact, assume
$s\geq4\log\bigl(5n+(4\pi)^4\bigr)$. Then the following hold:
\begin{enumerate}
\item Let $Z\subset X$ be an irreducible projective subvariety
of dimension $m>0$, and set $d_m=2\lfloor m/2\rfloor+1$.
Its middle intersection Betti number
$b_m^{\mathrm{IH}}(Z)=\dim_{\Qbb}IH^m(Z,\Qbb)$ satisfies
\[
  b_m^{\mathrm{IH}}(Z)>
  \frac1{d_m}\left(\frac{e^{s/16}}{4\pi}\right)^m.
\]
\item For $Z$ as in part~(i) and every nonzero semisimple complex
local system $\Lcal$ of rank $r$ on a smooth dense Zariski-open
subset of $Z$, one has
\[
  \frac{\dim_{\Cbb}IH^m(Z,\Lcal)}r>
  \frac1{d_m}\left(\frac{e^{s/16}}{4\pi}\right)^m.
\]
\end{enumerate}
\end{maincorollary}
\end{samepage}

\begin{introremark}
For smooth $Z$, part~(i) bounds its ordinary middle Betti number
$b_m(Z)$. If $\Lcal$ is defined on all of smooth $Z$, part~(ii)
bounds $\dim_{\Cbb}H^m(Z,\Lcal)/r$.
\end{introremark}

 A collection of projective subvarieties of
$X$ is called bounded if every member occurs as a fiber of a flat
projective family of subvarieties of $X$ over a scheme of finite
type over $\Cbb$. The inequalities used in the proof of Corollary~\ref{cor:betti} also give a bound on the degrees of projective subvarieties, yielding the following boundedness result:

\begin{maincorollary}[Bounded families of supports]
\label{cor:bounded-family}
Let $X=\Gamma\backslash\Ball^n$ be a finite-volume complex ball
quotient with neat lattice. Fix an integer $1\leq m\leq n$ and
$B>0$. Each of the following collections forms a bounded family
of subvarieties of $X:$
\begin{enumerate}
\item The irreducible projective $m$-dimensional subvarieties
$Z\subset X$ with $b_m^{\mathrm{IH}}(Z)\leq B$.
\item The irreducible projective $m$-dimensional subvarieties
$Z\subset X$ admitting a nonzero semisimple complex local system
$\Lcal$ of rank $r$ on a smooth dense Zariski-open subset such that
\[
  \frac{\dim_{\Cbb}IH^m(Z,\Lcal)}r\leq B.
\]
\end{enumerate}
\end{maincorollary}

 A cofinal normal tower of $X$ is a nested sequence of finite
\'etale Galois covers $X_i=\Gamma_i\backslash\Ball^n$, where
$\Gamma_i\triangleleft\Gamma$ are finite-index normal subgroups
satisfying $\Gamma_{i+1}\subseteq\Gamma_i$ and
$\bigcap_i\Gamma_i=\{1\}$. Along such a tower, the systoles
$\sys(X_i)$ tend to infinity
\cite[Proposition~3.21]{MemarianVolumes}. The lower bounds in
Theorem~\ref{thm:systolic} therefore tend to infinity as well,
giving the following consequence:

\begin{maincorollary}\label{cor:tower}
Let $X_i=\Gamma_i\backslash\Ball^n$ be a sequence
of finite covers of $X$. If
$\sys(X_i)\to\infty$, then
\[
\min\left\{
  \chi(X_i,\Pcal):
  \Pcal\in\Perv(X_i,\Cbb),\
  \dim\operatorname{Supp}\Pcal>0
\right\}
\xrightarrow{i\to\infty}\infty.
\]
In particular, this holds for every cofinal
normal tower of finite covers of $X$.
\end{maincorollary}

\newpage

\subsection*{Related work}
\pdfbookmark[2]{Related work}{related-work}

We discuss related work in two categories:
\begin{enumerate}
 \item \textbf{Euler characteristics of perverse sheaves.}
Franecki--Kapranov \cite{FraneckiKapranov} proved nonnegativity
of Euler characteristics of perverse sheaves on semiabelian
varieties.
Liu, Maxim, and Wang \cite[Corollary~3.4]{LiuMaximWang} proved the strict positivity
for nonzero perverse sheaves on proper ball quotients, while Arapura and Patel
\cite{ArapuraPatel} established nonnegativity of the orbifold
Euler characteristics of perverse sheaves on moduli stacks of
curves and principally polarized abelian varieties in
characteristic zero.
The first author and Tran \cite{AkdenizliTran} proved
nonnegativity on finite-volume Hermitian locally symmetric
spaces of noncompact type, and strict positivity for perverse
sheaves with full support.
More generally, Arapura--Mese--Patel \cite{ArapuraMesePatel}
proved strict positivity for perverse sheaves with full support
on smooth subvarieties with nonzero Euler
characteristic, provided the ambient variety admits a smooth
projective compactification with nef logarithmic cotangent
bundle.
Theorem~\ref{thm:criterion} provides strict positivity for every
nonzero perverse sheaf under a stronger positivity hypothesis on the
logarithmic cotangent bundle, without requiring full support.
Its application to ball quotients includes sheaves with singular
or nonproper support.

\item \textbf{Betti numbers.}
\pdfbookmark[3]{Betti numbers}{related-betti}
L\"uck's approximation theorem
\cite[Theorem~0.1]{LueckApproximation} identifies limits of Betti
numbers divided by covering degree along cofinal normal towers with
the von Neumann $L^2$-Betti numbers of the universal cover.
Yeung \cite{YeungBettiTower} obtained effective estimates for spaces
of holomorphic forms along cofinal normal towers of compact locally symmetric
domains. Di Cerbo--Stern \cite[Theorem~2]{DiCerboSternBetti}
give effective upper bounds for Betti numbers divided by volume
away from the middle degree on compact complex-hyperbolic manifolds.

These results describe the cohomology of the ambient manifolds.
Corollary~\ref{cor:betti} gives a lower bound for the middle
intersection Betti number of each irreducible projective
subvariety $Z\subset X$ in terms of the ambient systole alone.
The subvariety may be singular and may have any positive
dimension, and the ambient quotient may be noncompact.

In particular, let $X_i\to X$ be a cofinal normal tower of
finite covers, and let $Z_i\subset X_i$ be any sequence of
irreducible projective subvarieties of dimensions $m_i>0$.
Since $\sys(X_i)\to\infty$, Corollary~\ref{cor:betti} gives
\[
  b_{m_i}^{\mathrm{IH}}(Z_i)\longrightarrow\infty.
\]
When the $Z_i$ are smooth, their ordinary middle Betti
numbers therefore tend to infinity. 
\end{enumerate}

\subsection*{Strategy of proofs}
\pdfbookmark[2]{Strategy of proofs}{strategy-of-proofs}

Sections~\ref{sec:perverse} and~\ref{sec:positivity} establish the strict
positivity criterion. The logarithmic index theorem of Wu--Zhou
\cite{WuZhou} expresses $\chi(Y,\Pcal)$ as a sum of logarithmic
conormal intersection numbers with positive integral multiplicities.
We compare these numbers with the corresponding degrees after a negative
rational boundary twist of the logarithmic cotangent bundle. The
Chern-class twisting formula and nefness make the boundary contributions
nonnegative, while ampleness of the twist gives a strictly positive
term. This proves Theorem~\ref{thm:criterion}, and hence
Theorem~\ref{thm:positive} by \cite{Memarian}.

Sections~\ref{sec:boundary} and~\ref{sec:complete} give quantitative
bounds reflecting the geometry of the supports. For nonproper
subvarieties, whose closures meet the boundary, we keep a suitable term
in the twisting expansion to obtain a bound in terms of the boundary
degree. For projective subvarieties contained in $X$, a curvature
estimate shows that a negative canonical twist of the cotangent bundle
is nef. We then compare conormal intersections with canonical degrees
and use the subvariety-volume estimates of Hwang--To \cite{HwangTo}
and the second author \cite{MemarianVolumes}. Combining the two cases
proves Theorem~\ref{thm:systolic}.

\section*{Acknowledgment}
\pdfbookmark[1]{Acknowledgment}{acknowledgment}
The authors would like to thank Donu Arapura for several helpful discussions related to this work.
\paragraph{AI disclosure.}
GPT-5.6 Sol was used to assist
with language editing and polishing the text.
\section{Background on perverse sheaves}\label{sec:perverse}

Let $Y$ be a smooth complex algebraic variety of dimension $n$,
and let $\Ybar$ be a smooth projective compactification with simple
normal crossings boundary $D=\Ybar\setminus Y$. We recall the
definitions and properties of perverse sheaves used below, and the
logarithmic index formula that reduces our results to conormal
intersections.

Let $D^b_c(Y,\Cbb)$ be the bounded derived category of algebraically
constructible complexes of $\Cbb$-vector spaces. For
$\mathcal F\in D^b_c(Y,\Cbb)$, its support $\Supp\mathcal F$ is the
union of the supports of its cohomology sheaves, and its Euler
characteristic is
\[
  \chi(Y,\mathcal F)
  =\sum_k(-1)^k\dim_{\Cbb}\mathbb H^k(Y,\mathcal F).
\]
For an adapted Whitney stratification of $Y$, write $i_S:S\hookrightarrow Y$
for the inclusion of a stratum $S$. The complex $\mathcal F$ is
perverse, with middle perversity, if
\[
  \mathcal H^k(i_S^*\mathcal F)=0\quad(k>-\dim S),
  \qquad
  \mathcal H^k(i_S^!\mathcal F)=0\quad(k<-\dim S)
\]
for every stratum $S$. Here $\mathcal H^k$ denotes the $k$th
cohomology sheaf, and $i_S^*$ and $i_S^!$ are the usual and
exceptional inverse images. We denote this category by
$\Perv(Y,\Cbb)$; see \cite[\S2]{BBD}.

If $i:Z\hookrightarrow Y$ is a closed immersion with $Z$ smooth
of dimension $m$, and $\Lcal$ is a finite-rank complex local system
on $Z$, then
$i_*\Lcal[m]$ is perverse. In particular,
$\chi(Y,i_*\Cbb_Z[m])=(-1)^m\chi_{\mathrm{top}}(Z)$.
Subvarieties are understood to be reduced algebraic subsets of their
ambient varieties.

More generally, if $Z\subset Y$ is irreducible of dimension $m$
and $j:Z_{\mathrm{reg}}\hookrightarrow Z$ is its smooth locus,
then its intersection complex
\[
  IC_Z:=j_{!*}\Cbb_{Z_{\mathrm{reg}}}[m]
\]
is perverse. If $i:Z\hookrightarrow Y$ is the closed inclusion,
then $i_*IC_Z$ is perverse and
\[
  \chi(Y,i_*IC_Z)=(-1)^m\chi_{\mathrm{IH}}(Z).
\]
Here $IH^j(Z,\Cbb):=\mathbb H^{j-m}(Z,IC_Z)$ is ordinary
intersection cohomology, and
\[
 \chi_{\mathrm{IH}}(Z)
 =\sum_j(-1)^j\dim_{\Cbb}IH^j(Z,\Cbb).
\]

For an irreducible subvariety $Z\subset Y$, let $Z_{\rm reg}$
denote its smooth locus. The conormal bundle $T^*_{Z_{\rm reg}}Y$
consists of covectors at points of
$Z_{\rm reg}$ that vanish on its tangent spaces. Its closure
\[
  \Lambda_Z:=\overline{T^*_{Z_{\rm reg}}Y}^{\,T^*Y}
\]
is an $n$-dimensional conic Lagrangian subvariety of $T^*Y$.
Square brackets denote fundamental cycle classes. The
characteristic cycle of a perverse sheaf $\Pcal$ is effective
and can be written as
\begin{equation}\label{eq:CC}
  \CC(\Pcal)=\sum_\alpha\nu_\alpha[\Lambda_{Z_\alpha}],
  \qquad \nu_\alpha\in\mathbb Z_{>0},
\end{equation}
where the sum runs over its nonzero components and the $Z_\alpha$
are irreducible subvarieties contained in $\Supp\Pcal$.
Every irreducible component $Z$ of $\Supp\Pcal$ contributes its
conormal with multiplicity equal to the positive rank of the
shifted local system on a dense smooth open subset of $Z$.
See \cite[Corollary~5.2.24]{Dimca} and
\cite[Proposition~3.2]{ArapuraWang}; for characteristic cycles,
see \cite[\S4.3]{Dimca} and \cite{KashiwaraIndex}.

Set $E=\Omega^1_{\Ybar}(\log D)$, so that $E|_Y=T^*Y$.
For an irreducible subvariety $Z\subset Y$, let
\(
  \Lambda_Z^{\log}:=\overline{\Lambda_Z}^{\,E}
\)
be the closure of its conormal variety in the total space of $E$.
This is an $n$-dimensional conic subvariety. Write
$0_E(\Ybar)\hookrightarrow{} E$ for the zero section and define
\begin{equation}\label{eq:Ilog}
  \Ilog(Z):=[\Lambda_Z^{\log}]\cdot[0_E(\Ybar)]\in\mathbb Z.
\end{equation}
Here the dot denotes the degree of the zero-section
intersection, a zero-cycle on $\Ybar$.

The logarithmic closure of the characteristic cycle in
\eqref{eq:CC}, with its multiplicities, is
\[
  \overline{\CC(\Pcal)}^{\,E}
   =\sum_\alpha\nu_\alpha[\Lambda_{Z_\alpha}^{\log}].
\]
Wu--Zhou's logarithmic Dubson--Kashiwara index theorem
\cite[Theorem~1.6]{WuZhou} states precisely that
\begin{equation}\label{eq:index}
  \chi(Y,\Pcal)
  =[\overline{\CC(\Pcal)}^{\,E}]\cdot[0_E(\Ybar)]
  =\sum_\alpha\nu_\alpha\Ilog(Z_\alpha).
\end{equation}
 When $D=\varnothing$, this
is the usual Dubson--Kashiwara index formula
\cite[Theorem~4.3.25]{Dimca}. The positive
multiplicities reduce our results to bounds for $\Ilog(Z)$.

\papersection{Chern class inequalities and strict positivity}\label{sec:positivity}

\subsection{Chern class inequalities for rational twists}

We recall the Chern-class inequalities needed for conormal
intersections and then prove the strict positivity criterion. The main
inputs are the positivity results of Lazarsfeld
\cite[Chapter~8]{LazarsfeldPositivityII} and Fulton--Lazarsfeld
\cite[\S1 and Theorem~2.3]{FultonLazarsfeld}.

We use the quotient convention: $\PP(G)$ parametrizes
one-dimensional quotients of $G$, with tautological quotient line
bundle $\Ocal_{\PP(G)}(1)$. For a vector bundle $G$ on a projective
variety $\Ybar$, write $p_G:\PP(G)\to\Ybar$ for the projection. Given a rational Cartier
class $\delta$, denote the $\Qbb$-twisted bundle by
$G\langle\delta\rangle$. We write $c_k$ for the $k$th Chern class,
with $c_0=1$. The twist is nef or ample when the rational class
$c_1(\Ocal_{\PP(G)}(1))+p_G^*\delta$ is nef or ample.
If $G$ has rank
$r$, shifting its formal Chern roots gives
\begin{equation}\label{eq:twist}
  c_k(G\langle\delta\rangle)
   =\sum_{i=0}^k\binom{r-i}{k-i}c_i(G)\delta^{k-i};
  \qquad
  c_r(G\langle\delta\rangle)
   =\sum_{j=0}^r c_{r-j}(G)\delta^j.
\end{equation}

For an integral subvariety $V\subset\Ybar$, $[V]\in A_{\dim V}(\Ybar)$
denotes its fundamental cycle class in the Chow group; an effective
cycle is a nonnegative integral combination of such classes. An
integral sign denotes the degree of the resulting zero-cycle, with
rational coefficients when a $\Qbb$-twist is involved.

\begin{lemma}\label{lem:mixed}
If $G\langle\delta\rangle$ is nef on a projective variety $\Ybar$,
$\alpha$ is an effective $t$-cycle, and $H$ is a nef rational Cartier
class, then
\[
  \int_\alpha c_k(G\langle\delta\rangle)
      H^{t-k}\geq0
  \qquad(0\leq k\leq\min\{t,\operatorname{rk}G\}).
\]

\end{lemma}
\begin{proof}
By linearity it suffices to consider $\alpha=[V]$ with $V$ integral of
dimension $t$. Restrict all classes to $V$. Assume first that
$H$ is ample. Choose $N>0$ such that $NH$ is very ample,
and choose $t-k$ general divisors in $|NH|$. They meet
properly and give an effective $k$-cycle
$N^{t-k}H^{t-k}\cap[V]
=\sum_\nu m_\nu[W_\nu]$, where $m_\nu>0$ and
$\dim W_\nu=k$. The projection formula gives
\[
 N^{t-k}\int_Vc_k(G\langle\delta\rangle)H^{t-k}
 =\sum_\nu m_\nu\int_{W_\nu}
      c_k(G|_{W_\nu}\langle\delta|_{W_\nu}\rangle)\geq0.
\]
Every restricted bundle is nef, and its Chern degree is
nonnegative by \cite[Theorem~8.2.1]{LazarsfeldPositivityII}.
For a singular $W_\nu$, apply that theorem on a resolution.
The pulled-back twist remains nef, and the projection formula
preserves its Chern degree.
This includes $t=k$, with an empty divisor intersection. For nef
$H$, fix an ample class $A$ and define the polynomial
\[
 f(\varepsilon)=\int_Vc_k(G\langle\delta\rangle)
     (H+\varepsilon A)^{t-k}.
\]
For every rational $\varepsilon>0$, the class $H+\varepsilon A$
is ample, so $f(\varepsilon)\geq0$ by the case already proved.
Letting $\varepsilon\to0^+$ gives
$f(0)=\int_Vc_k(G\langle\delta\rangle)H^{t-k}\geq0$.

\end{proof}

A line in a fiber of $G$ is a quotient of
the dual fiber $G^\vee$. Consequently a cone $C\subset G$ has a
projectivization inside $\PP(G^\vee)$. Write
$\pi:\PP(G^\vee)\to\Ybar$ for the bundle projection.
Dualizing the tautological
quotient on $\PP(G^\vee)$ gives the exact sequence
\begin{equation}\label{eq:universal}
  0\longrightarrow\Ocal_{\PP(G^\vee)}(-1)
  \longrightarrow\pi^*G\longrightarrow Q_G\longrightarrow0.
\end{equation}

\begin{lemma}\label{lem:universal}
Let $G\langle\delta\rangle$ have rank $r>0$ on a projective variety
$\Ybar$, and let $\alpha$ be a nonzero effective $(r-1)$-cycle on
$\PP(G^\vee)$. If $G\langle\delta\rangle$ is nef, then
\[
  \int_\alpha c_{r-1}
       (Q_G\langle\pi^*\delta\rangle)\geq0.
\]
Moreover, the inequality is strict if $G\langle\delta\rangle$ is ample.
\end{lemma}

\begin{proof}
For $r=1$, the integral is the degree of a nonzero
effective zero-cycle, hence positive. Assume $r>1$. If $G\langle\delta\rangle$ is nef, then
$Q_G\langle\pi^*\delta\rangle$ is nef as a quotient
of $\pi^*G\langle\pi^*\delta\rangle$.
The inequality follows from Lemma~\ref{lem:mixed}.

Suppose $G\langle\delta\rangle$ is ample.
Choose an ample divisor $H$ on $\Ybar$ and a small
rational $\varepsilon>0$ such that
$G\langle\delta-\varepsilon H\rangle$ is nef.
Then \(R=Q_G\langle\pi^*(\delta-\varepsilon H)\rangle\)
is nef. For any integral component $V$ of $\alpha$,
the twist gives
\[
  \int_V c_{r-1}(Q_G\langle\pi^*\delta\rangle)
  =\sum_{k=0}^{r-1}\varepsilon^{r-1-k}
    \int_V c_k(R)(\pi^*H)^{r-1-k}.
\]
Every summand is nonnegative by Lemma~\ref{lem:mixed}.

Set $W=\pi(V)$, $t=\dim W$, and $q=r-1-t$.
The generic fiber $V_\eta\subset\PP^{r-1}$ has
dimension $q$. By the universal sequence
\eqref{eq:universal}, the restriction of $c_q(R)$
to this fiber is $h^q$, where
$h=c_1(\Ocal_{\PP^{r-1}}(1))$. Thus
\[
  \pi_*(c_q(R)\cap[V])=e[W],
  \qquad e=\int_{V_\eta}h^q>0.
\]
The summand with $k=q$ is therefore
\[
  \varepsilon^t\int_V c_q(R)(\pi^*H)^t
  =\varepsilon^t e\,H^t\cdot[W]>0.
\]
Hence the degree is positive on every component
of $\alpha$, proving the strict inequality.
\end{proof}

\subsection{A criterion for strict positivity}

We use the notation of Section~\ref{sec:perverse}. In particular,
$Y=\Ybar\setminus D$ and $E=\Omega^1_{\Ybar}(\log D)$.
The following proposition is the main ingredient in
Theorem~\ref{thm:criterion}.

\begin{proposition}
\label{prop:weighted-cone}
Let $\Ybar$ be smooth projective of dimension $n>0$, and let
$D=\sum_{i=1}^rD_i$ be a reduced SNC divisor. Suppose that
$E|_{D_i}$ is nef for each $i$ and that
\[
 E\langle-A\rangle\quad\text{is ample},\qquad
 A=\sum_{i=1}^r a_iD_i,\qquad a_i\in\Qbb_{>0}.
\]
Then $\Ilog(Z)>0$ for every irreducible subvariety
$Z\subset Y=\Ybar\setminus D$. If $D$ is empty, the assertion
holds under the assumption that $E$ is ample.
\end{proposition}

\begin{proof}
The residue map for each component $D_i$ is a surjection
$E|_{D_i}\twoheadrightarrow\Ocal_{D_i}$. After twisting by
$-A$, this becomes a quotient of rationally twisted bundles
\[
 E\langle-A\rangle|_{D_i}
 \twoheadrightarrow\Ocal_{D_i}\langle-A|_{D_i}\rangle.
\]
A quotient of an ample rationally twisted vector bundle is
ample. Hence $L_i:=\Ocal_{D_i}(-A)$ is an ample
$\Qbb$-line bundle.

If $Z$ is a point, then $\Ilog(Z)=1$. Suppose
$0<m:=\dim Z<n$. Put $B=\PP(\Lambda_Z^{\log})
\subset\PP(E^\vee)$, let $\pi:\PP(E^\vee)\to\Ybar$,
and write $Q=Q_E$ for the universal quotient bundle.
By formula
\cite[\S1, (1.5)]{FultonLazarsfeld},
\begin{equation}\label{eq:cone-formula}
 \Ilog(Z)=\int_B c_{n-1}(Q).
\end{equation}
Set
\[
 I_A(Z):=\int_B c_{n-1}(Q\langle-\pi^*A\rangle).
\]
Lemma~\ref{lem:universal} applied to the ample rational twist
$E\langle-A\rangle$ gives $I_A(Z)>0$.

For each $i$, let $H_i=\pi^{-1}(D_i)$ and
$\beta_i=H_i\cdot[B]$. Since $B$ is the closure of a
cycle lying over $Y$, no component of $B$ is contained in
$H_i$, and $\beta_i$ is effective. The
twisting formula and
$A^k=\sum_i a_i D_i A^{k-1}$ give
\begin{equation}\label{eq:weighted-comparison}
\begin{aligned}
 \Ilog(Z)-I_A(Z)
 &=\sum_{i=1}^r a_i\sum_{k=1}^{n-1}
     \int_{\beta_i}
     c_{n-1-k}(Q|_{H_i})
     (\pi|_{H_i})^*c_1(L_i)^{k-1}
 \;\geq 0.
\end{aligned}
\end{equation}
Indeed, $Q|_{H_i}$ is nef as a quotient of
$(\pi|_{H_i})^*(E|_{D_i})$, and
$(\pi|_{H_i})^*c_1(L_i)$ is nef on $H_i$.
Every summand is nonnegative by Lemma~\ref{lem:mixed}.
Therefore $\Ilog(Z)\geq I_A(Z)>0$.

Finally, when $Z=Y$, the logarithmic conormal is the zero
section, so $\Ilog(Y)=\int_{\Ybar}c_n(E)$. Positivity for ample
rational twists \cite[Theorem~8.2.4]{LazarsfeldPositivityII} gives
$\int_{\Ybar}c_n(E\langle-A\rangle)>0$.
The same expansion gives
\begin{equation}\label{eq:weighted-zero-section}
\begin{aligned}
 \int_{\Ybar}c_n(E)-\int_{\Ybar}c_n(E\langle-A\rangle)
 &=\sum_{i=1}^r a_i\sum_{k=1}^{n}
   \int_{D_i}c_{n-k}(E|_{D_i})c_1(L_i)^{k-1}
 \\&\geq0.
\end{aligned}
\end{equation}
Here the inequality follows from nefness of $E|_{D_i}$. This proves the assertion
for all $Z$. When $D=\varnothing$, the sums are empty,
and the same strict positivity assertions apply directly.
\end{proof}

\begin{proof}[Proof of Theorem~\ref{thm:criterion}]
Under the hypotheses of the theorem, the restrictions
$E|_{D_i}$ are nef. Proposition~\ref{prop:weighted-cone}
shows that $\Ilog(Z)>0$ for every irreducible
$Z\subset Y$. A nonzero perverse sheaf has a nonzero
effective characteristic cycle. The logarithmic index
formula~\eqref{eq:index} therefore gives
\[
 \chi(Y,\Pcal)=\sum_\alpha\nu_\alpha\Ilog(Z_\alpha)>0.
\]
\end{proof}

\begin{proof}[Proof of Theorem~\ref{thm:positive}]
Let $\Xbar$ be the smooth toroidal compactification
of the noncompact ball quotient $X$, and let
$D=\Xbar\setminus X$. The logarithmic cotangent bundle
$E=\Omega^1_{\Xbar}(\log D)$ is nef by
\cite[Corollary~G]{Memarian}, while
$E\langle-aD\rangle$ is ample for any rational
$0<a<d/(4\pi)$ by \cite[Theorem~A]{Memarian}, where
$d>0$ denotes the uniform cusp depth. Taking
$a_i=a$ for every irreducible component of $D$,
Theorem~\ref{thm:criterion} applies.
\end{proof}

\section{Inequalities for supports meeting the boundary}\label{sec:boundary}
Let $X=\Gamma\backslash\Ball^n$ be a finite-volume ball quotient
with neat lattice, $\Xbar$ its smooth projective toroidal
compactification, and $D=\Xbar\setminus X$ its boundary. When $X$
is compact, we take $\Xbar=X$ and $D=\varnothing$. Write
$E=\Omega^1_{\Xbar}(\log D)$, and let $\overline Z\subset\Xbar$
denote the closure of an irreducible subvariety $Z\subset X$.
We refine the weighted comparison of
Proposition~\ref{prop:weighted-cone} to obtain the boundary-degree
bound used in Theorem~\ref{thm:systolic}.

For $0<m:=\dim Z<n$, write
$P_Z=\PP(\Lambda_Z^{\log})\subset\PP(E^\vee)$ and
$\pi:P_Z\to\overline Z$. The cone formula~\eqref{eq:cone-formula}
applies with $B=P_Z$ and $Q=Q_E$. Over the generic point of $Z$,
the fiber of $\pi$ is $\PP^{n-m-1}$.

For $D\neq\varnothing$, let $d$ denote the uniform cusp depth in the
normalization of \cite[Definition~2.3]{Memarian}. The boundary conormal
line bundle $\Ocal_D(-D)$ is ample by the residue-quotient
argument of Proposition~\ref{prop:weighted-cone}, applied to
$E\langle-aD\rangle$ (see also \cite[\S1]{MokProjective}).
For a positive-dimensional $Z$, define its boundary degree by
\begin{equation}\label{eq:bD}
  b_D(Z):=\deg_{\Ocal_D(-D)}\bigl(D\cdot[\overline Z]\bigr)
  =D\cdot(-D)^{m-1}\cdot[\overline Z].
\end{equation}
Here $D\cdot[\overline Z]$ is the effective boundary cycle on $D$.
If $\overline Z$ meets $D$, this cycle is nonzero, so ampleness gives
$b_D(Z)\in\mathbb Z_{>0}$; otherwise $b_D(Z)=0$.
\begin{lemma}\label{lem:boundary-degree}
Suppose $X$ is noncompact, and let $Z\subset X$ be a
positive-dimensional irreducible subvariety. Write
$m=\dim Z$. For every rational number $a$ with $0<a<d/(4\pi)$,
one has
\begin{equation}\label{eq:boundary-component-bound}
  \Ilog(Z)>a^m b_D(Z).
\end{equation}
\end{lemma}
\begin{proof}
For the stated values of $a$, the positivity hypotheses of
Proposition~\ref{prop:weighted-cone} hold with $A=aD$, as recalled
in the proof of Theorem~\ref{thm:positive}.
First suppose $m<n$, and put $q=n-m-1$, $Q=Q_E$ and $H=\pi^*D$.
On a general $\PP^q$-fiber of $P_Z\to\overline Z$, the universal
sequence~\eqref{eq:universal} gives $\int_{\PP^q}c_q(Q)=1$.
Consequently,
\begin{equation}\label{eq:conormal-pushforward}
 \pi_*\bigl(c_q(Q)\cap[P_Z]\bigr)=[\overline Z].
\end{equation}
Specializing~\eqref{eq:weighted-comparison} to $A=aD$, the twisted
cone degree is strictly positive and all terms in the comparison
are nonnegative. Retaining the term with $k=m$ and using the
projection formula, we obtain
\[
 \Ilog(Z)>a^m\int_{H\cdot[P_Z]}
 c_q(Q|_H)c_1(\Ocal_H(-H))^{m-1}
 =a^m b_D(Z).
\]
If $m=n$, then $Z=X$. The zero-section comparison
\eqref{eq:weighted-zero-section}, retaining the term $k=n$,
similarly gives
\[
 \Ilog(X)>a^n\int_D c_1(\Ocal_D(-D))^{n-1}
 =a^n b_D(X).
\]
\end{proof}

\begin{corollary}\label{cor:boundary}
Let $Z\subset X$ be an irreducible subvariety of dimension $m>0$
whose closure meets $D$. Then
\[
  \Ilog(Z)\geq
  \left\lceil\left(\frac d{4\pi}\right)^m b_D(Z)\right\rceil
  \geq1.
\]
\end{corollary}
\begin{proof}
Lemma~\ref{lem:boundary-degree} gives
$\Ilog(Z)>a^m b_D(Z)$ for every rational $0<a<d/(4\pi)$.
Letting rational $a$ tend to $d/(4\pi)$ gives
$\Ilog(Z)\geq(d/(4\pi))^m b_D(Z)$.
Since $\Ilog(Z)$ is an integer by~\eqref{eq:Ilog}, the claimed
ceiling follows.
\end{proof}

\section{Systolic bounds}\label{sec:complete}

Let $X$, its smooth projective toroidal compactification $\Xbar$,
and $D$ be as in Section~\ref{sec:boundary}. For an irreducible
subvariety $Z\subset X$, write $\overline Z$ for its closure in
$\Xbar$. We first estimate the conormal intersection for projective
$Z$. Proposition~\ref{prop:complete} compares this intersection with
the canonical degree and then applies the volume estimates in terms
of $\sys(X)$. We combine this with the estimate for
$\overline Z\cap D\neq\varnothing$ to obtain the bound for every
$Z\subset X$ and prove the main results.

\subsection{Bisectional curvature}
\label{subsec:curvature-conventions}

In the ball coordinates $z=(z_1,\ldots,z_n)$, with
$|z|^2=\sum_j|z_j|^2$, use the Bergman metric normalized to have
holomorphic sectional curvature $-1$. Its K\"ahler form and
Hermitian coefficients are
\[
 \omega_{\Ball^n}
 =-2i\partial\bar\partial\log(1-|z|^2)
 =i\sum_{j,k}g_{j\bar k}\,dz_j\wedge d\bar z_k,
 \qquad
 g_{j\bar k}
 =2\left(\frac{\delta_{jk}}{1-|z|^2}
       +\frac{\bar z_jz_k}{(1-|z|^2)^2}\right).
\]
Here $\delta_{jk}$ is the Kronecker delta, and $g$ is the
Hermitian metric on the holomorphic tangent bundle. It descends
to $X$, with K\"ahler form $\omega_X$. The associated Riemannian
metric is
\[
 ds^2=2\sum_{j,k}g_{j\bar k}\,dz_j\,d\bar z_k.
\]
Its restriction to a complex geodesic disk is
$4|d\zeta|^2/(1-|\zeta|^2)^2$, of Gaussian curvature $-1$.

Let $h=g^\vee$ be the induced metric on the cotangent bundle.
In the holomorphic coframe $(dz_1,\ldots,dz_n)$, its coefficient
matrix is $(g^{\mathsf t})^{-1}$, so
\[
 h_{k\bar l}
 =\frac12(1-|z|^2)(\delta_{kl}-z_k\bar z_l).
\]
We use the Chern curvature tensor with components
\[
 (R_h)_{i\bar j k\bar l}
 =-\partial_i\partial_{\bar j}h_{k\bar l}
  +\sum_{a,b}h^{a\bar b}
       (\partial_i h_{k\bar b})
       (\partial_{\bar j}h_{a\bar l}),
\]
where $\sum_b h_{k\bar b}h^{a\bar b}=\delta_k^a$.
For nonzero $v\in T_{X,x}^{1,0}$ and
$\xi\in\Omega_{X,x}^1$, the holomorphic bisectional curvature
of this metric on the cotangent bundle is defined to be
\[
 B_h(v,\xi)
 =\frac{R_h(v,\bar v;\xi,\bar\xi)}
        {\|v\|_g^2\|\xi\|_h^2}.
\]
The following inequality is well known for the unit ball and 
we include a computation for completeness:

\begin{lemma}\label{lem:bisectional-curvature}
For the metric above and nonzero
$v\in T_{X,x}^{1,0}$ and $\xi\in\Omega_{X,x}^1$,
\begin{equation}\label{eq:bisectional-bound}
 B_h(v,\xi)
 =\frac12\left(1+
   \frac{|\xi(v)|^2}{\|v\|_g^2\|\xi\|_h^2}\right)
 \geq\frac12.
\end{equation}
\end{lemma}
\begin{proof}
By homogeneity under holomorphic isometries, it suffices to
compute at the origin. The metric on the cotangent bundle has
Taylor expansion
\[
 h_{k\bar l}(z)
 =\frac12\delta_{kl}
  -\frac12\delta_{kl}|z|^2
  -\frac12 z_k\bar z_l+O(|z|^4).
\]
At the origin, $h_{k\bar l}(0)=\frac12\delta_{kl}$ and all first
derivatives vanish at the origin. The Chern curvature formula gives
\begin{equation}\label{eq:ball-curvature-tensor}
 \begin{aligned}
 (R_h)_{i\bar j k\bar l}(0)
 &=-\partial_i\partial_{\bar j}h_{k\bar l}(0)\\
 &=\frac12(\delta_{ij}\delta_{kl}
                +\delta_{ik}\delta_{jl}).
 \end{aligned}
\end{equation}
Since $g(0)=2I$ and $h(0)=\frac12 I$, contraction yields
\[
 R_h(v,\bar v;\xi,\bar\xi)
 =\frac12\left(
   \|v\|_g^2\|\xi\|_h^2+|\xi(v)|^2\right).
\]
This identity is invariant under holomorphic isometries, so it
holds everywhere on the ball and descends to $X$. Dividing by
$\|v\|_g^2\|\xi\|_h^2$ proves
\eqref{eq:bisectional-bound}.
\end{proof}

\subsection{The curvature twist and canonical degree}
\label{subsec:curvature-twist}

Equip $K_X=\det\Omega_X^1$ with the determinant metric $h_K$.
For a Hermitian line bundle $(M,h_M)$ and a nonvanishing local
holomorphic frame $e$, its first Chern form is
$c_1(M,h_M)=-\frac{i}{2\pi}\partial\bar\partial\log\|e\|_{h_M}^2$.
Taking the trace of the curvature formula above gives
\begin{equation}\label{eq:canonical-class}
 R_{h_K}(v,\bar v)=\frac{n+1}{2}\|v\|_g^2,
 \qquad
 c_1(K_X,h_K)=\frac{n+1}{4\pi}\omega_X,
\end{equation}
where $R_{h_K}$ denotes the Chern curvature of $h_K$.

\begin{lemma}\label{lem:curvature}
For every projective subvariety $Z\subset X$,
the $\Qbb$-twisted bundle
\[
 \left.\Omega_X^1
 \left\langle-\frac{K_X}{n+1}\right\rangle\right|_Z
\]
is nef.
\end{lemma}
\begin{proof}
By Lemma~\ref{lem:bisectional-curvature} and
\eqref{eq:canonical-class}, the curvature of the twisted metric
satisfies
\[
 R_h(v,\bar v;\xi,\bar\xi)
 -\frac{R_{h_K}(v,\bar v)}{n+1}\|\xi\|_h^2
 =\frac12|\xi(v)|^2\geq0.
\]
The twist has nonnegative bisectional curvature. The induced
tautological rational line bundle on $\PP(\Omega_X^1|_Z)$ has
semipositive curvature, hence nonnegative degree on every curve.
This proves nefness.
\end{proof}

\subsection{Canonical degree and the systolic estimate}
\label{subsec:complete-systole}
All lengths use $ds^2$. Write $\operatorname{dist}_g$ for its
distance on the ball. For $\gamma\in\Gamma$, its translation
length is
$\ell(\gamma)=\inf_{z\in\Ball^n}\operatorname{dist}_g(z,\gamma z)$.
An element is hyperbolic if it preserves a geodesic axis and
translates along it by a positive distance. The systole of $X$ is
defined to be
\begin{equation}\label{eq:systole}
 \sys(X)=\inf_{\substack{\gamma\in\Gamma\\\gamma\text{ hyperbolic}}}
            \ell(\gamma).
\end{equation}
For a finite-volume quotient, this infimum is positive and attained
\cite[Proposition~3.9 and Remark~3.19]{MemarianVolumes}.
Parabolic elements do not enter the definition. For $x\in X$
and a lift $\tilde x\in\Ball^n$, the injectivity radius is
\[
 \operatorname{inj}_x(X)=\frac12\inf_{\gamma\in\Gamma\setminus\{1\}}
       \operatorname{dist}_g(\tilde x,\gamma\tilde x),
\]
independently of the lift. When $X$ is proper, this gives
$\operatorname{inj}_x(X)\geq\sys(X)/2$.
The induced volume of $Z$ is
$\operatorname{vol}_X(Z)=\int_{Z_{\mathrm{reg}}}\omega_X^m/m!$.

\begin{proposition}\label{prop:complete}
Let $Z\subset X$ be an irreducible projective subvariety,
and write $m=\dim Z>0$. Then
\begin{equation}\label{eq:complete-bound}
 \Ilog(Z)\geq
 \sinh^{2m}\!\left(\frac{\sys(X)}8\right).
\end{equation}
\end{proposition}

\begin{proof}
We first show that
\begin{equation}\label{eq:proper-comparison}
  \Ilog(Z)\geq\frac{K_X^m\cdot[Z]}{(n+1)^m}
  =\frac{m!}{(4\pi)^m}\operatorname{vol}_X(Z).
\end{equation}

Set
\(F=\Omega_X^1|_Z\), and \(\lambda=\frac{c_1(K_X|_Z)}{n+1}.\)
By \eqref{eq:canonical-class}, $\lambda$ is nef,
since it is represented by the positive form
$\omega_X/(4\pi)$. Moreover,
Lemma~\ref{lem:curvature} shows that
$F\langle-\lambda\rangle$ is nef.

Suppose $m<n$. Let $\pi:P_Z\to Z$ be the projectivized
conormal variety and $Q_F$ the universal quotient
from \eqref{eq:universal}. Then
$R=Q_F\langle-\pi^*\lambda\rangle$ is nef.
The cone formula \eqref{eq:cone-formula} and the
twist formula \eqref{eq:twist} give
\[
  \Ilog(Z)
  =\sum_{k=0}^{n-1}
    \int_{P_Z}c_k(R)(\pi^*\lambda)^{n-1-k}.
\]
Every summand is nonnegative by Lemma~\ref{lem:mixed}.
For $q=n-m-1$, the general fiber of $\pi$ is
$\PP^q$, on which $c_q(R)$ has degree one. Hence
\(
  \pi_*\bigl(c_q(R)\cap[P_Z]\bigr)=[Z],
\)
and keeping the term $k=q$ gives
\[
  \Ilog(Z)\geq
  \int_{P_Z}c_q(R)(\pi^*\lambda)^m
  =\lambda^m\cdot[Z].
\]
If $m=n$, then $Z=X$, and the same argument gives
\[
  \Ilog(X)=
  \sum_{k=0}^n\int_X
  c_k(F\langle-\lambda\rangle)\lambda^{n-k}
  \geq\lambda^n.
\]
This proves \eqref{eq:proper-comparison}, with the
volume identity following from \eqref{eq:canonical-class}.

If $X$ is proper, Hwang--To \cite[Theorem~1.1]{HwangTo} and
$\operatorname{inj}_x(X)\geq\sys(X)/2$ give
$\frac{m!}{(4\pi)^m}\operatorname{vol}_X(Z)
\geq\sinh^{2m}(\sys(X)/4)$.
If $X$ is nonproper, the volume estimate
\cite[Theorem~5.2]{MemarianVolumes} implies, in particular,
the weaker bound
$\frac{m!}{(4\pi)^m}\operatorname{vol}_X(Z)
\geq\sinh^{2m}(\sys(X)/8)$.
Combining these estimates with \eqref{eq:proper-comparison}
proves the proposition.
\end{proof}

\subsection{Proofs of the main results}\label{sec:consequences}

\begin{lemma}\label{lem:depth-systole}
Suppose $X$ is noncompact, write $s=\sys(X)$, and assume
$s\geq4\log\bigl(5n+(4\pi)^4\bigr)$. Then
\begin{equation}\label{eq:depth-systole}
  d>e^{s/16}>4\pi.
\end{equation}
\end{lemma}
\begin{proof}
This is \cite[Corollary~3.15]{MemarianVolumes}.
\end{proof}

\begin{proposition}\label{prop:systolic-conormal}
Let $Z\subset X$ be an irreducible subvariety of dimension $m>0$,
and write $s=\sys(X)$. If $X$ is noncompact, assume
$s\geq4\log\bigl(5n+(4\pi)^4\bigr)$. Then
\begin{equation}\label{eq:systolic-conormal}
  \Ilog(Z)>\left(\frac{e^{s/16}}{4\pi}\right)^m.
\end{equation}
\end{proposition}
\begin{proof}
If $\overline Z$ meets $D$, Corollary~\ref{cor:boundary},
Lemma~\ref{lem:depth-systole}, and $b_D(Z)\geq1$ give the result.
If $X$ is noncompact and $Z$ is projective, use
Proposition~\ref{prop:complete}. The hypothesis gives
$s>16\log(4\pi)$, and
\[
  \sinh^2(s/8)\geq\frac1{16}e^{s/4}
  >\frac{e^{s/16}}{4\pi}.
\]
For compact $X$, put $t=e^{s/16}$. If $t<4\pi$, use the positive
integer $\Ilog(Z)\geq1$ from Proposition~\ref{prop:weighted-cone} and integrality.
If $t\geq4\pi$, the compact estimate in the proof of
Proposition~\ref{prop:complete} gives
$\Ilog(Z)\geq\sinh^{2m}(s/4)$, which proves the claim.
\end{proof}

\begin{proof}[Proof of Theorem~\ref{thm:systolic}]
Let $Z\subset X$ be the stated support component of $\Pcal$.
The logarithmic index formula and Proposition~\ref{prop:weighted-cone}
give $\chi(X,\Pcal)\geq\Ilog(Z)$.
For compact $X$, the compact estimate in the proof of
Proposition~\ref{prop:complete} proves part~(i).
For noncompact $X$, Proposition~\ref{prop:systolic-conormal}
proves part~(ii).
\end{proof}

\begin{proof}[Proof of Corollary~\ref{cor:tower}]
Let $m_i>0$ be the dimension of any positive-dimensional
irreducible support component of $\Pcal_i$.
If $X$ is nonproper, Theorem~\ref{thm:systolic}(ii) gives,
for all sufficiently large $i$,
\[
 \chi(X_i,\Pcal_i)>
 \left(\frac{e^{\sys(X_i)/16}}{4\pi}\right)^{m_i}
 \geq\frac{e^{\sys(X_i)/16}}{4\pi}.
\]
If $X$ is proper, Theorem~\ref{thm:systolic}(i) gives
\[
 \chi(X_i,\Pcal_i)\geq
 \sinh^{2m_i}\!\left(\frac{\sys(X_i)}4\right)
 \geq\sinh^2\!\left(\frac{\sys(X_i)}4\right)
\]
for all sufficiently large $i$. Both right-hand sides tend to
infinity as $\sys(X_i)\to\infty$.
For a cofinal normal tower, this systolic divergence follows from
\cite[Proposition~3.21]{MemarianVolumes}.
\end{proof}

\begin{corollary}\label{cor:length}
Let $X$ be as in Theorem~\ref{thm:positive}.
Then every $\Pcal\in\Perv(X,\Cbb)$ satisfies
\(
\chi(X,\Pcal)\geq\operatorname{length}(\Pcal),
\)
where $\operatorname{length}(\Pcal)$ denotes its
Jordan--Hölder length.
\end{corollary}

\begin{proof}
Let $0=\Pcal_0\subset\Pcal_1\subset\cdots
\subset\Pcal_r=\Pcal$ be a Jordan--Hölder filtration.
By additivity of Euler characteristics and
Theorem~\ref{thm:positive},
\(
\chi(X,\Pcal)
=\sum_{i=1}^r
\chi(X,\Pcal_i/\Pcal_{i-1})
\geq r.
\)
\end{proof}

\begin{remark}
By applying Theorem~\ref{thm:systolic} to the simple
constituents of a Jordan--Hölder filtration and using
additivity of Euler characteristics, one obtains
corresponding lower bounds scaled by the number of
positive-dimensional constituents.
\end{remark}

\section{Middle Betti numbers and boundedness}
\label{sec:betti}
In this section we prove Corollary~\ref{cor:betti} and Corollary~\ref{cor:bounded-family}. 
\begin{proof}[Proof of Corollary~\ref{cor:betti}]
Let $Z\subset X$ be an irreducible projective subvariety of
dimension $m>0$, and set $d_m=2\lfloor m/2\rfloor+1$.
Let $\Lcal$ be a nonzero semisimple complex local system of rank
$r$ on a smooth dense Zariski-open subset $j:U\hookrightarrow Z$.
Write
\[
  IC_Z(\Lcal)=j_{!*}(\Lcal[m]),
  \qquad IH^q(Z,\Lcal)=\mathbb H^{q-m}(Z,IC_Z(\Lcal)).
\]
For the inclusion $i:Z\hookrightarrow X$, the perverse sheaf
$\Pcal=i_*IC_Z(\Lcal)$ has effective characteristic cycle, with
the conormal of $Z$ occurring with multiplicity $r$.
The logarithmic index formula~\eqref{eq:index} and
Proposition~\ref{prop:weighted-cone} give
\[
  \chi(Z,IC_Z(\Lcal))=\chi(X,\Pcal)\geq r\Ilog(Z)>0.
\]

Put $b_q=\dim_{\Cbb}IH^q(Z,\Lcal)$, and let $\eta$ be an ample
class on $Z$. Hard Lefschetz for semisimple perverse sheaves
\cite[Theorem~2.1.2 and Corollary~2.2.1]{deCataldoSemisimple} gives
\begin{equation}\label{eq:IH-HL}
  \eta^k:IH^{m-k}(Z,\Lcal)
  \xrightarrow{\ \sim\ }IH^{m+k}(Z,\Lcal)
  \qquad(k\geq0).
\end{equation}
For $0\leq2a\leq m$, the composition
\[
  IH^{m-2a}(Z,\Lcal)\xrightarrow{\eta^a}IH^m(Z,\Lcal)
  \xrightarrow{\eta^a}IH^{m+2a}(Z,\Lcal)
\]
is the isomorphism $\eta^{2a}$. Hence
$b_{m-2a}=b_{m+2a}\leq b_m$. There are $d_m$ indices with the
same parity as $m$ between $0$ and $2m$, so
\[
  \chi(Z,IC_Z(\Lcal))
  =\sum_{q=0}^{2m}(-1)^{q-m}b_q
  \leq\sum_{\substack{0\leq q\leq2m\\q\equiv m\pmod{2}}}b_q
  \leq d_m b_m.
\]
For every finite-volume ball quotient with neat lattice, we have
therefore obtained the comparison
\begin{equation}\label{eq:IH-conormal}
  \frac{\dim_{\Cbb}IH^m(Z,\Lcal)}r
  \geq\frac{\chi(Z,IC_Z(\Lcal))}{r d_m}
  \geq\frac{\Ilog(Z)}{d_m}>0.
\end{equation}
This comparison holds without a systolic hypothesis.
Under the assumptions of the corollary,
Proposition~\ref{prop:systolic-conormal} now gives part~(ii).
Taking constant coefficients gives part~(i).
\end{proof}

\begin{proof}[Proof of Corollary~\ref{cor:bounded-family}]
Choose a small rational $\varepsilon>0$ and an integer $\ell>0$
such that
\[
  A=\ell\bigl(K_{\Xbar}+(1-\varepsilon)D\bigr)
\]
is a very ample Cartier divisor, using Bakker--Tsimerman
\cite[Theorem~A]{BakkerTsimermanKodaira}. For compact $X$, take
$A=\ell K_X$.
Let $Z\subset X$ be a projective support of dimension $m$, set
$d_m=2\lfloor m/2\rfloor+1$, and let $\Lcal$ be a semisimple
local system of rank $r$ as in part~(b). One has
$A|_Z=\ell K_X|_Z$. The comparison~\eqref{eq:proper-comparison}
and \eqref{eq:IH-conormal} give
\begin{equation}\label{eq:IH-degree-bound}
  A^m\cdot[Z]
  \leq\bigl(\ell(n+1)\bigr)^m\Ilog(Z)
  \leq d_m\bigl(\ell(n+1)\bigr)^m
  \frac{\dim_{\Cbb}IH^m(Z,\Lcal)}r.
\end{equation}
Consequently, $\deg_A Z\leq d_m\bigl(\ell(n+1)\bigr)^mB$
in both collections, using constant coefficients for part~(a).

Since $A$ is very ample, the degree bound places all such
$Z$ among integral subvarieties of a fixed projective space
of dimension $m$ and bounded degree. By the standard
boundedness theorem, only finitely many Hilbert polynomials
occur. Therefore, these subvarieties are parametrized by a finite
union $T$ of Hilbert schemes of $\Xbar$, with a universal
flat projective family \cite{GrothendieckHilbert}.

Let $T^\circ\subset T$ be the locus parametrizing
geometrically integral fibers disjoint from $D$.
This locus is open, since geometric integrality is
open in flat projective families, and the condition
of avoiding $D$ is open by properness.
Restricting the universal family to $T^\circ$ gives
the required flat projective family in
$X\times T^\circ$, with $T^\circ$ of finite type.
\end{proof}

\begin{samepage}
\bigskip
\begingroup\small
\hypersetup{urlcolor=black}
\noindent\textsc{Mert Akdenizli}\par\nobreak
\noindent Department of Mathematics, Purdue University,
West Lafayette, IN 47907, USA\par\nobreak
\noindent\textit{Email address}: \href{mailto:makdeniz@purdue.edu}{\nolinkurl{makdeniz@purdue.edu}}\par
\medskip
\noindent\textsc{Soheil Memariansorkhabi}\par\nobreak
\noindent Department of Mathematics, Purdue University,
West Lafayette, IN 47907, USA\par\nobreak
\noindent\textit{Email address}: \href{mailto:smemaria@purdue.edu}{\nolinkurl{smemaria@purdue.edu}}\par
\endgroup
\end{samepage}
\end{document}